\documentclass[11pt]{amsart}
\usepackage{graphicx}
\usepackage{epsfig}
\usepackage{amsmath}
\usepackage{amssymb}
\usepackage{amsfonts}
\usepackage{cite, tocvsec2}
\usepackage{amsfonts, amstext, amsmath, amssymb, mathrsfs, amsthm}
\usepackage{amsopn,amscd,xspace,enumerate,multicol}

\usepackage{xspace}
\usepackage[breaklinks=true]{hyperref}

\newcommand\lie[1]{\mathfrak{#1}}

\newcommand{\g}{\lie{g}}

\newcommand{\s}{\lie{sl}}
\newcommand{\C}{\mathbb{C}}

\newcommand{\Q}{\mathbb{Q}}
\newcommand{\Z}{\mathbb{Z}}

\newcommand\op[1]{{\rm{#1}}} 

\newtheorem{theorem}{Theorem}

\newtheorem{corollary}[theorem]{Corollary}
\newtheorem{proposition}{Proposition}

\begin{document} 

\title[The center of $u_q(\s_3)$, $l=5$]{The center of the small quantum group $u_q(\s_3)$ for $l=5$: a numerical computation.} 
\author{Steven Glenn Jackson,  Anna Lachowska} 
\maketitle 

\section{Introduction} 

The small quantum group is a finite dimensional Hopf algebra $u_q(\g)$ associated to a semisimple Lie 
algebra $\g$ and a root of unity $q^l =1$. It was defined by \cite{Lus} as an analog of the restricted enveloping algebra $u_p(\g)$ 
over a field of positive characteristic. 
Despite the fact that $u_q(\g)$ is a finite dimensional algebra over $\C$ that has been 
studied extensively for over 20 years, 
the structure and even the dimension of its center remains unknown, except in the case of 
$\g = \s_2$. In the latter case, the answer was first found in \cite{Ker}: the dimension of the center 
of $u_q(\s_2)$ with $q$ a primitive root of unity of odd degree $l \geq 3$, equals $\frac{3l-1}{2}$, which 
is unexpectedly large (the number of inequivalent irreducible representations of $u_q(\s_2)$ is $l$).  

For higher rank, there are several results on the structure of the regular block  of the center 
of $u_q(\g)$, whose structure for large enough $l$ is independent of $l$ by \cite{AJS}. 
The main theorem in  \cite{BL} gives a complete description of the Hochschild cohomology of the regular block of $u_q(\g)$ 
for a semisimple $\g$ in terms of the sheaf cohomology of certain polyvector fields over the Springer 
resolution. This description does not easily translate into the combinatorial terms. However it allows us to conclude that
 the regular block of the center 
of $u_q(\g)$ contains a subalgebra of 
dimension $2|W|-1$, where $W$ is the Weyl group associated with $\g$. The existence of this subalgebra
was predicted in \cite{La}, and for $\g = \s_2$ it is known 
 to coincide with the whole center of the regular block.  Further, in 
  \cite{LQ} the dimension of 
the center of the regular block of $u_q(\s_3)$ was computed, which led to a conjecture for the structure of the regular block of the center in type A. Still, there is no answer for the dimension of the entire center of $u_q(g)$ beyond $g=\s_2$.

In the present note we describe the result of a numerical computation of the dimension of the center 
of $u_q(\s_3)$ at a particular root of unity, $l=5$, performed in 2011. This dimension is $57$. Now 
in view of the result in \cite{LQ} for the regular block of the center, we can make sense of this number. Namely, 
using the linkage principle \cite{APW},  
 we obtain the block 
decomposition of the algebra $u_q(\s_3)$. It contains blocks of three kinds: regular 
blocks, singular parabolic blocks and the Steinberg block. 
The number of each kind of blocks is determined by the parameter $l$.  The Steinberg block 
is a matrix algebra and its center is one-dimensional. The dimension of the regular block of 
the center is $16$ by \cite{LQ}. Therefore, we can find the dimension of the center of the singular parabolic 
block for $\g=\s_3$ (it is $6$), and derive the block decomposition of the center for $u_q(\s_3)$. 

The paper is organized as follows. 
Section 2 contains the definition of the small quantum group $u_q(\s_3)$ with mild restrictions on $l$, and a PBW basis of this algebra over the cyclotomic field. In Section 3 we present our strategy and the result of the numerical computation 
of the dimension of the center of this algebra for $l=5$. In Section 4 we interpret this result in terms of the block decomposition 
of the center of $u_q(\s_3)$.

{\bf Acknowledgements}  The authors are grateful to Bryan Ford and DEDIS group at the Computer Science Department of Yale University for allowing us to use their computational facilities. 

\section{The small quantum $\s_3$} 

We follow \cite{Lus} for the definition of the small quantum group $u_q = u_q(\s_3)$ corresponding to the Lie algebra $\s_3$ and 
an $l$-th root of unity, where $l >3 $ is an odd integer prime to $3$.

Let $U_v(\s_3)$ be the Drinfeld-Jimbo quantum group over $\Q(v)$ corresponding to the Lie algebra $\s_3$. 
Let $a_{ij}$ be the $2\times 2$ Cartan matrix with $a_{11} = a_{22} =2$, $a_{12} = a_{21} =-1$. 
Then $U_v(\s_3)$  is defined by the generators  $E_i, F_i, K_i, K_i^{-1}$, where $i= 1,2$ and the relations 
\begin{equation}  \label{U_v} 
\begin{array}{ll} 
K_i K_j = K_j K_i   & K_i K_i^{-1} = K_i^{-1} K_i =1  \\[2mm] 
K_i E_j = v^{a_{ij}} E_j K_i    & K_i F_j = v^{-a_{ij}} F_j K_i   \\
E_i F_j - F_j E_i = \delta_{ij} \frac{K_i - K_i^{-1}}{v-v^{-1}}  &  \\[2mm] 
E_i^2 E_j  - (v-v^{-1}) E_i E_jE_i + E_j E_i^2 =0 ,  &   i \neq j \\[2mm] 
F_i^2 F_j  - (v-v^{-1}) F_i F_jF_i + F_j F_i^2 =0 ,  &   i \neq j \\ 
\end{array}  
\end{equation} 
 
Let $A = \Z[v, v^{-1}]$ and let $U^{\op{div}}$ denote the $A$-subalgebra of $U_v(\s_3)$ generated by 
\[ E_i^{(n)} = \frac{E_i^n}{[n]!} , \quad F_i^{(n)}  = \frac{F_i^n}{[n]!} , \quad K_i, K_i^{-1}, \quad i=1,2 \] 
where $[n] = \frac{v^n - v^{-n}}{v-v^{-1}}$, and $[n]! = [1] \cdot [2] \cdot \ldots [n]$. 
Let $q \in \C$ be a primitive $l$-th root of unity, and let $U_q = U^{\op{div}} \otimes_A \Z[q, q^{-1}]$.   
Then $(F_i^{(1)})^l = (E_i^{(1)})^l =0$ and $K_i^{2l} =1$ in $U_q$ for $i=1,2$. Let $u_q$ be the $\Z[q, q^{-1}]$-subalgebra 
of $U_q$ generated by $E_i^{(1)}, F_i^{(1)}, K_i, K_i^{-1}$, for $i=1,2$, modulo the two-sided ideal generated 
by the central elements $K_1^l-1, K_2^l-1$. It is a finite dimensional Hopf algebra. 

Following \cite{Lus}, we introduce the elements $E_3 = q^{-1} E_1 E_2 - E_2 E_1$ and $F_3 = q F_2 F_1 - F_1 F_2$ that correspond to the highest root of $\s_3$. 
Then the algebra $u_q(\s_3)$ admits the following basis over $\Q(q)$: 
\begin{equation} \label{PBW}
 \left\{ \prod F_1^{n_1} F_3^{n_3} F_2^{n_2}  \prod K_1^{k_1} K_2^{k_2} \prod E_2^{m_2}  E_3^{m_3} E_1^{m_1}  \right\}_{0 \leq n_i, m_i, k_j \leq l-1, \; i=1,2,3, \; j=1,2} . \end{equation} 
 Thus the dimension of the algebra $u_q$ over $\Q(q)$ is $l^8$. 

\section{The numerical computation of the dimension of the center} 

Let $l=5$, so that $q$ is the $5$th root of unity, $q \neq 1$. 
We are interested in the dimension of the center $z(u_q)$. 

To search for the central elements in the basis (\ref{PBW}) it will be convenient to reformulate 
the commutation relations (\ref{U_v}) 
using the generators $E_3, F_3$. We obtain in addition to (\ref{U_v}):  
\medskip
\begin{equation} 
\begin{array}{ll} 
E_1 E_2 = q E_2 E_1 + q E_3  &    E_2 E_1 = q^{-1} E_1 E_2 - E_3 \\[2mm] 
E_1 F_3 = F_3 E_1 - F_2 K_1^{-1}  &  E_2 F_3 = F_3 E_2 + q F_1 K_2   \\[2mm]
 E_1 E_3 = q^{-1} E_3 E_1    &  E_2 E_3 = q E_3 E_2 \\[2mm] 
 F_1F_2 = q F_2 F_1 - F_3   &  F_2 F_1 = q^{-1} F_1 F_2 + q^{_1} F_3  \\[2mm] 
 F_1 E_3 = E_3 F_1 + K_1 E_2  &  F_2 E_3 = E_3 F_2 - q^{-1} K_2^{-1} E_1 \\[2mm] 
 F_1 F_3 = q^{-1} F_3 F_1  &   F_2 F_3 = q F_3 F_2 \\[2mm] 
 K_1 E_3  = q E_3 K_1 & K_2 E_3  = q E_3 K_2  \\[2mm]
 K_1 F_3 = q^{-1} F_3 K_1  &  K_2 F_3 = q^{-1} F_3 K_2 \\[2mm] 
 E_3 F_3  = F_3 E_3  + \frac{K_1 K_2 - K_1^{-1} K_2^{-1}}{q - q^{-1}}  &   \\[2mm]  
 \end{array} 
 \end{equation} 
 
We briefly describe the computation strategy. 
The $K_1, K_2$-invariant subspace $Z_{K_1, K_2}= (u_q)^{(0,0)} $ of $u_q$ is spanned by the elements of weight $(0,0)$ of the basis (\ref{PBW}). Equivalently, these elements satisfy the conditions
\[ 2n_1 -n_2 + n_3 = 2m_1 - m_2 + m_3,  \quad \quad 2n_2 - n_1 + n_3 = 2m_2 - m_1 +m_3 .\] 
The dimension of this space is $8125$. Let $\{v_s\}_{s=1}^{8125}$ be an ordered set of monomials from (\ref{PBW})
 that form a basis in $Z_{K_1, K_2}$. 
 
 Next we find the dimensions of the centralizer subspaces $Z_{E_1}, Z_{E_2}, Z_{F_1}, Z_{F_2} \subset Z_{K_1, K_2}$, 
 spanned by the elements in $Z_{K_1, K_2}$ that commute with $E_1, E_2, F_1, F_2$ respectively. 
 The elements $E_i v_s - v_s E_i$ have the weight $(2, -1)$ for $E_1$ and $(-1,2)$ for $E_2$, while the elements 
 $F_i v_s - v_s F_i$ have the weight $(-2, 1)$ for $F_1$ and $(1,-2)$ for $F_2$. It is easy to compute that the dimension of each of the subspaces $(u_q)^{(2,-1)}$,  $(u_q)^{(-1,2)}$, $(u_q)^{(-2,1)}$, $(u_q)^{(1,-2)}$,  of $u_q$ is $7500$. In each of the subspaces $(u_q)^\nu$,  where $\nu = (2,-1), (-1,2), (-2,1), (1,-2)$,  choose an ordered  set of monomials $ \{u^\nu_j \}_{j=1}^{7500}$ from (\ref{PBW}) that form a basis in $(u_q)^\nu$.  

Then we can write  
\[ \begin{array}{ll} 
 E_1 v_s - v_s E_1 = \sum_{j=1}^{7500} M^{E_1}_{js} \, u^{(2,-1)}_j   & \;\;\;\; s = 1, \ldots 8125 \\[2mm] 
 E_2 v_s - v_s E_2 = \sum_{j=1}^{7500} M^{E_2}_{js} \, u^{(-1,2)}_j   & \;\;\;\; s = 1, \ldots 8125 \\[2mm] 
F_1 v_s - v_s F_1 = \sum_{j=1}^{7500} M^{F_1}_{js}  \, u^{(-2,1)}_j   & \;\;\;\; s = 1, \ldots 8125 \\[2mm] 
F_2 v_s - v_s F_2 = \sum_{j=1}^{7500} M^{F_2}_{js}  \, u^{(1,-2)}_j   & \;\;\;\;  s = 1, \ldots 8125 \\[2mm] 
\end{array} \]  
where $M^{E_1}$, $M^{E_2}$, $M^{F_1}$, and $M^{F_2}$  are matrices of size $8125 \times 7500$ with entries in $\Q(q)$. 
Now the problem is reduced to  
\begin{enumerate} 
\item Computing the matrices $M^{E_1}$, $M^{E_2}$, $M^{F_1}$, and $M^{F_2}$; 
\item Finding the null spaces  $Z_{E_i} = \op{Ker} \,M^{E_i}$, $Z_{F_i} = \op{Ker}\,M^{F_i}$ of each matrix as subspaces 
in $Z_{K_1, K_2}$. 
\item Finding the intersection of these subspaces.
\end{enumerate}   
We will write  $Z_{E_1, F_1} = Z_{E_1} \cap Z_{F_1}$, etc. 

The computation in SAGE produced the following results: \\
(Here we can include more on the computation, the packages used, particular tricks, etc.) 
\[ \begin{array}{l} 
\op{dim} Z_{E_1} = \op{dim} Z_{E_2} = \op{dim} Z_{F_1} = \op{dim} Z_{F_2} = 3325 \\[2mm] 
\op{dim} Z_{E_1, E_2} = \op{dim} Z_{F_1, F_2} = 195 \\[2mm] 
\op{dim} Z_{E_1, F_1} = \op{dim} Z_{E_2, F_2} =  925 \\[2mm] 
\op{dim} Z_{E_1, F_2} = \op{dim} Z_{E_2, F_1} = 721 \\[2mm] 
\op{dim} Z_{E_1, E_2, F_1} = \op{dim} Z_{E_1, E_2, F_2} = \op{dim} Z_{F_1, F_2, E_1} = \op{dim} Z_{F_2, E_1, E_2} = 81 \\[2mm] 
\op{dim} \, z(u_q) = \op{dim} Z_{E_1, E_2, F_1, F_2} = 57 .
\end{array} \]

 \section{The block decomposition of $z(u_q)$.} 
 
 Let $l$ be an odd integer prime to $3$ and such that $l >3$. As a finite dimensional associative algebra, $u_q = u_q(\s_3)$ decomposes as a direct sum of two-sided ideals, or blocks. This 
 decomposition is determined by a version of the linkage principle stated in \cite{APW, APW2}. Let $\alpha_1, \alpha_2$ be the simple roots of $\s_3$, $Q$ the root lattice and $P$ the weight lattice. Recall (\cite{Lus}) that the inequivalent 
 simple modules $\{L(\lambda)\}$ over $u_q$ are parametrized by the highest weights $\lambda \in P$ such that 
 \[ 0  \leq  \langle \lambda , \alpha_i \rangle \leq l-1 , \quad \quad i=1,2 .\] 
 We will denote this set by $P/lP$. 
 
 \begin{proposition} \cite{APW, APW2}.  If two simple $u_q$-modules $L(\lambda)$ and $L(\mu)$ appear in the filtration of an 
 indecomposable $u_q$-module $M$, then $\lambda$ and $\mu$ belong to the same orbit of the shifted action of 
 the extended affine Weyl group $\hat{W}_{l,P}$ :  
 \[  \mu \in  \hat{W}_{l, P} \cdot \lambda , \quad \quad \hat{W}_{l, P} = W \ltimes lP, \] 
 where $w \cdot \lambda = w(\lambda + \rho) - \rho$. 
 \end{proposition}  
 
 \begin{corollary} \cite{La} The block decomposition of $u_q$ is determined by the orbits of the shifted action of $\hat{W}_{l,P}$ in 
 $P/lP$. Let $\bar{\chi}$ denote the set of representatives of the orbits of $\hat{W}_{l,P}$ in $P/lP$. Then we have 
 \[  u_q  = \oplus_{\nu \in \bar{\chi}} \;  (u_q)_\nu . \] 
 \end{corollary} 
 As a left $u_q$-module, each block $(u_q)_\nu$ is a direct sum of indecomposable projective modules with composition 
 factors of highest weights in the set $\hat{W}_{l,P} \cdot \nu \in  P/lP$.  
 
 For $u_q$, the set $P/lP$ is a parallelogram generated by the $(l-1)$-multiples of the fundamental weights, $P/lP = \{ (k_1, k_2\} : 0\leq k_1, k_2 \leq l-1$. With respect to the $\hat{W}_{l,P}\cdot$-action there are three kinds of weights: the Steinberg weight $\mu = (l-1) \rho$ that constitutes in itself an $\hat{W}_{l,P}\cdot$-orbit; $3(l-1)$ the weights of the form $\{  (k_1, k_2) : k_1 + k_2 + 2 =l\} \cup \{(i, l-1) : 0 \leq i \leq l-2\} \cup \{(l-1,j) : 0 \leq j \leq l-2\}$ that are stabilized by one of the reflections in $\hat{W}_{l,P}$; and the remaining $(l-1)(l-2)$ weights in $P/lP$, that have a trivial stabilizer in $\hat{W}_{l,P}$. 
 
 \medskip
Therefore the orbits in $\bar{\chi}$ can be divided into three groups: 
 \medskip
 \begin{center} 
 \begin{tabular}{|c|c|c|c|}  \hline
       $ \mu \in \bar{\chi}$  &  stabilizer &  order  of the orbit  & number of orbits \\[2mm]  \hline 
       regular                &   $\{1\}$        &                     6             & $\frac{(l-1)(l-2)}{6}$      \\[2mm] \hline  
       singular   parabolic  &  proper subgroup &            3          &   $ l-1 $    \\[2mm] \hline 
       $\mu = (l-1)\rho$      &   $\hat{W}_{l,P}$    &          1             &   1                     \\[2mm] \hline  
 \end{tabular}       
 \end{center} 
 \bigskip 
 
 The center of $u_q(\s_3)$ decomposes according to the same principle: 
 \[ z(u_q) = \oplus_{\nu \in \bar{\chi}} \;  (z(u_q))_\nu .  \] 
 The blocks of the center corresponding to the orbits of the same kind (regular or singular parabolic) are isomorphic. 
 Therefore we obtain for the dimension of the center
 \[ \op{dim}\, z(u_q) = \frac{(l-1)(l-2)}{6}  \cdot  \op{dim}\, z_q^{\op{reg}} + (l-1) \cdot  \op{dim} \, z_q^{\op{par}} + 1 ,\] 
 where $z_q^{\op{reg}}$ and $z_q^{\op{par}}$ denote respectively the center of the regular and parabolic blocks. 
By   Theorem 3.6 in \cite{LQ}, the dimension of $z_q^{\op{reg}}$ is $16$. Therefore, we have for $l=5$:
 \[ 57 =  2 \cdot 16 + 4 \cdot  \op{dim} \, z_q^{\op{par}} + 1 ,\]
 which implies that in this case 
 \[ \op{dim}\, z_q^{\op{par}} = 6. \] 
 
 Following the philosophy in \cite{AJS} we can expect that the singular parabolic block of the center of $u_q$ will have 
 the same dimension for all odd $l >3$ prime to $3$. Then we have that the dimension of the center for $u_q = u_q(\s_3)$ 
 should be equal to 
 \[ \op{dim}\, z(u_q) = \frac{(l-1)(l-2)}{6} \cdot 16 + (l-1) \cdot 6 + 1  .\]

\end{document}